\documentclass{amsart}
\usepackage{amssymb,latexsym}
\usepackage{amsfonts}

\newcommand{\ba}{\begin{array}}
\newcommand{\ea}{\end{array}}

\def \qed{\cqfd}

\def\qed{\vbox{\hrule
\hbox{\vrule\hbox to 5pt{\vbox to 8pt{\vfil}\hfil}\vrule}\hrule}}

\newcommand{\beg}{\begin{eqnarray*}}
\newcommand{\begn}{\begin{eqnarray}}
\newcommand{\en}{\end{eqnarray*}}
\newcommand{\enn}{\end{eqnarray}}

\begin{document}
\title{Existence of approximate Hermitian-Einstein structures on semi-stable Higgs bundles}
\author{Jiayu Li}
\address{School of Mathematical Sciences,\\
University of Science and Technology of China,\\
Hefei, 230026,\\ and AMSS, CAS, Beijing, 100080, P.R. China\\} \email{jiayuli@ustc.edu.cn}
\author{Xi Zhang}
\address{School of Mathematical Sciences,\\
University of Science and Technology of China,\\
Hefei, 230026,P.R. China\\ } \email{mathzx@ustc.edu.cn}
\thanks{The authors were supported in part by NSF in
China, No.11071212, No.11131007, No.10831008 and No. 11071236.}

\begin{abstract} In this paper, using  Donaldson's heat flow,
 we show that  the semi-stability of a Higgs bundle over a compact K\"ahler manifold implies the existence of approximate Hermitian-Einstein structure on the Higgs bundle.
\end{abstract}

\maketitle

\section{Introduction}
\setcounter{equation}{0}

\hspace{0.4cm}

Let $(M, \omega )$ be a compact K\"ahler manifold,  and  $E$ be a holomorphic vector bundle
over $M$.  The stability of holomorphic vector bundles, in the sense of Mumford-Takemoto,  was a well established
concept in algebraic geometry. A holomorphic vector bundle $E$ is
called  stable (semi-stable), if for every coherent sub-sheaf
$E'\hookrightarrow E$ of lower rank, it holds:
\begin{eqnarray}
\mu (E')=\frac{deg (E')}{rank E'}< (\leq ) \mu (E)=\frac{deg (E)}{rank E},
\end{eqnarray}
where $\mu (E')$ is called the slope of $E'$. In early 1980¡¯s, S. Kobayashi introduced the Hermitian-Einstein condition for holomorphic
bundles on Kahler manifolds.
A Hermitian metric $H$ in $E$ is said to be Hermitian-Einstein, if the curvature $F_{H}$ of the Chern connection $D_{H}$
 satisfies  the Einstein condition:
\begin{eqnarray}
\sqrt{-1}\Lambda_{\omega }F_{H}=\lambda Id_{E},
\end{eqnarray}
where $\Lambda_{\omega }$ denotes the contraction of differential
forms by K\"ahler form $\omega $, and the real constant $\lambda $ is
given by $\lambda  =\frac{2\pi}{Vol(M)} \mu (E)$.

The so-called Hitchin-Kobayashi correspondence asserts that holomorphic vector bundles over compact K\"ahler manifolds are polystable if and only if they admit a Hermitian-Einstein metric. This correspondence starts by Narasimhan and Seshadri (\cite{NS}) in the case of compact Riemannian surface. Kobayashi (\cite{Ko1}) and L\"ubke (\cite{Lu}) proved that a holomorphic bundle admits a Hermitian-Einstein metric must be polystable. The inverse problem was solved by Donaldson (\cite{D1}, \cite{D2}) for algebraic manifolds, by Uhlenbeck and Yau (\cite{UY}) for general K\"ahler manifolds.  Donaldson-Uhlenbeck-Yau theorem states that the stability of a holomorphic vector  bundle implies the existence of  Hermitian-Einstein metric. The classical Hitchin-Kobayashi
correspondence has several interesting and important
generalizations and extensions where some extra structures are
added to the holomorphic bundles, see references: \cite{Hi}, \cite{Si},\cite{Br1},\cite{GP},
\cite{BG}, \cite{AG1},
\cite{Bi}, \cite{BT},  \cite{LN1}, \cite{LN2}, \cite{LY}, \cite{M}, \cite{Ta}.

\medskip

 A holomorphic bundle $(E , \overline{\partial }_{E})$ coupled with one Higgs field $\theta \in \Omega^{1,0}(End(E))$
 which satisfying $\overline{\partial}_{E}\theta =0$ and $\theta \wedge \theta =0$
 will be called by a  Higgs bundle.
A Higgs bundle $(E, \theta )$ is called  Stable (Semi-stable) if the usual stability
condition $\mu (E' ) <\mu (E)$ ($\leq $) hold for all proper
$\theta$-invariant sub-sheaves.
A Hermitian metric $H$ in Higgs bundle $(E, \theta )$ is said to be
 Hermitian-Einstein  if the
curvature $F$ of the Hitchin-Simpson connection
 $D_{H, \theta }= D_{H} + \theta +\theta ^{\ast H}$  satisfies the Einstein condition, i.e
\begin{eqnarray}
\sqrt{-1}\Lambda_{\omega} (F_{H} +[\theta , \theta^{\ast H}])
=\lambda Id_{E},
\end{eqnarray}
where $F_{H}$ is the curvature of  Chern connection $D_{H}$, $\theta ^{\ast H}$ is the adjoint of $\theta $ with respect to the metric $H$.

Higgs bundles  first emerged twenty
years ago in  Hitchin's study of the self-duality
equations on a Riemann surface and in  Simpson¡¯s
 work on nonabelian
Hodge theory. Higgs bundles have a rich structure and
play a role in many different areas including gauge theory,
K\"ahler and hyperk\"ahler geometry, group representations and
nonabelian Hodge theory. In \cite{Hi} and \cite{Si}, it is proved  that a Higgs bundle admits the
Hermitian-Einstein metric iff it's Higgs poly-stable. This is  a
Higgs bundle version of the Donaldson-Uhlenbeck-Yau theorem.

We say a holomorphic vector bundle $E$ admits an approximate Hermitian-Einstein structure if for every positive $\epsilon
$, there is a  Hermitian metric $H$ such that
\begin{eqnarray}
\max _{M} |\sqrt{-1}\Lambda_{\omega }F_{H}-\lambda Id_{E}|_{H}<\epsilon .
\end{eqnarray}
Kobayashi (\cite{Ko2}) introduced  the notion of approximate Hermitian-Einstein structure in a holomorphic vector bundle, and he proved that a holomorphic vector bundle with an approximate Hermitian-Einstein structure must be  semi-stable. Furthermore, over projective algebraic manifolds, Kobayashi solved  the inverse problem, i.e. the semi-stability implies admitting an approximate Hermitian-Einstein structure, and he also conjectured that the result should be true for general compact K\"ahler manifolds.

In this article, we consider the existence problem of approximate Hermitian-Einstein structure in Higgs bundles. We will show that  the semi-stability of Higgs bundle implies the existence of approximate Hermitian-Einstein structures. In fact, we prove the following theorem.

\medskip

{\bf Theorem 1. } {\it If $(E, \theta )$ is a semi-stable Higgs bundle on K\"ahler manifold $(M, \omega )$, then it admits an approximate Hermitian-Einstein structure, i.e. for any $\epsilon >0$ there exists a Hermitian metric such that
\begin{eqnarray}
\max _{M} |\sqrt{-1}\Lambda_{\omega }(F_{H}+[\theta , \theta ^{\ast H}]) -\lambda Id|_{H} <\epsilon .
\end{eqnarray}
 }

\medskip

We will use the heat flow method to prove theorem 1. Simpson (\cite{Si}) introduced the Donaldson's heat flow to Higgs bundle case, i.e.  the following heat flow for Hermitian
metrics on the Higgs bundle $(E, \theta)$ with initial
metric $H_{0}$:
\begin{eqnarray}\label{D1}
H^{-1}\frac{\partial H}{\partial
t}=-2(\sqrt{-1}\Lambda_{\omega}(F_{H}+[\theta , \theta ^{\ast H}
])-\lambda Id_{E}).
\end{eqnarray}
Simpson proved the long time existence and uniqueness of solution for the above non-linear heat equation. By Uhlenbeck and Yau's result, that $L_{1}^{2}$ weakly holomorphic sub-bundle defines a coherent sub-sheaf, Simpson obtain an uniform $C^{0}$-estimate  of the long time solution of the above heat flow (\ref{D1}) when the Higgs bundle is stable, and show that the solution must convergence to a Hermtian-Einstein metric. In this article, we will follow Simpson's discussion to prove that, along the heat flow, the term $\max _{M} |\sqrt{-1}\Lambda_{\omega }(F_{H}+[\theta , \theta ^{\ast H}]) -\lambda Id|_{H}$ must convergence to zero under the assumption that the Higgs bundle is semi-stable. The correspondence between semistability and the existence of approximate
Hermitian-Einstein structure in the holomorphic vector bundle case has been studied recently by Jacob (\cite{Ja})
 using a technique developed by Buchdahl (\cite{Bu}) for the regularization of
sheaves in the case of compact complex surfaces. It should be point out that our discussion is different from that in \cite{Ja}. Recently,  Cardona (\cite{Ca}) obtain the result of theorem 1 in Rieamnn surface case by using a Donaldson functional approach analogous to that of Kobayashi \cite{Ko2}.

In \cite{BO}, Bruzzo and Otero proved that a Higgs bundle admitting  an approximate Hermitian-Einstein structure must be semi-stable. Combining Bruzzo and Otero's result and theorem 1, we know that, in Higgs bundles,  admitting an approximate Hermitian-Einstein structure and the semi-stability are equivalent. It is easy to check that if two Higgs bundles admit approximate Hermitian-Einstein structure, so does their tensor product; furthermore, if they are with the same slope, so does their Whitney sum. So, we have the following corollary.

\medskip

{\bf Corollary 2.} {\it
Let $(E_{1}, \theta _{1} )$ and $(E_{2}, \theta _{2} )$ be two semi-stable Higgs bundle on K\"ahler manifold $(M, \omega )$, then

\medskip

(1), $(E_{1}\otimes E_{2}, \theta )$ is a semi-stable Higgs bundle, where $\theta =\theta_{1}\otimes Id_{2} +Id_{1}\otimes \theta _{2}$;

(2),  If $\frac{deg (E_{1})}{rank E_{1}}=\frac{deg (E_{2})}{rank E_{2}}$, then $(E_{1}\oplus E_{2}, \theta )$ is also a semi-stable Higgs bundle, where $\theta =pr_{1}^{\ast }\theta_{1}+ pr_{2}^{\ast} \theta _{2}$ and $pr_{i}: E_{1}\oplus E_{2} \rightarrow E_{i}$ denote the natural projections.

}

\medskip

Another application of theorem 1 is the following  Bogomolov type inequality for semi-stable Higgs bundle. By Chern-Weil theory, we have
\begin{eqnarray}
\begin{array}{lll}
&&4\pi^{2}\int_{M} (2C_{2}(E)-\frac{r-1}{r}C_{1}(E)\wedge C_{1}(E))\frac{\omega^{n-2}}{(n-2)!}\\
&=&\int_{M}tr (F_{H, \theta }^{\bot}\wedge F_{H, \theta }^{\bot})\wedge \frac{\omega^{n-2}}{(n-2)!}\\
&=&\int_{M}|F_{H, \theta }^{\bot }|_{H}^{2}-|\Lambda_{\omega} F_{H, \theta }^{\bot}|_{H}^{2} \frac{\omega^{n}}{n!}\\
&\geq & -\int_{M}|\sqrt{-1}\Lambda_{\omega} F_{H, \theta }-\lambda Id -\frac{1}{r}tr (\sqrt{-1}\Lambda_{\omega} F_{H, \theta }-\lambda Id)Id|_{H}^{2} \frac{\omega^{n}}{n!},\\
\end{array}
\end{eqnarray}
where $F_{H, \theta }^{\bot}$ is the trace free part of $F_{H, \theta }$. If the Higgs bundle $(E, \theta )$ admits an approximate Hermitian-Einstein structure, we can choose a sequence of metric $H_{i}$ so that the last term of the above inequality convergenes to zero, so we obtain the following corollary.

\medskip

{\bf Corollary 3.} {\it
If $(E, \theta )$ is a semi-stable Higgs bundle on K\"ahler manifold $(M, \omega )$, then we
have the following Bogomolov type inequality
\begin{eqnarray}
\int_{M} (2C_{2}(E)-\frac{r-1}{r}C_{1}(E)\wedge C_{1}(E))\frac{\omega^{n-2}}{(n-2)!}\geq 0 .
\end{eqnarray}}

This paper is organized as follows. In Section 2, we recall some basic estimates for the Donaldson's heat flow in Higgs bundle. In section 3, we prove theorem 1.

\hspace{0.3cm}

\section{Analytic preliminaries and basic estimates }
\setcounter{equation}{0}

Suppose $H(t)$ is a solution of the above Donaldson's heat flow (\ref{D1}) with initial metric $K$, and let
$h(t)=K^{-1}H(t)$, then (\ref{D1}) can be written as
\begin{eqnarray}
\frac{\partial h}{\partial
t}=-2\sqrt{-1}h\Lambda_{\omega}(F_{K}+\overline{\partial
}_{E}(h^{-1}\partial_{K}h)+[\theta, h^{-1}
\theta^{\ast K} h])+2\lambda h.
\end{eqnarray}

Furthermore, by an appropriate conformal change, we can assume that the initial metric $K$ satisfies
\begin{eqnarray}
tr (\sqrt{-1}\Lambda_{\omega}(F_{K}+[\theta , \theta^{\ast K}
])-\lambda Id_{E})=0.
\end{eqnarray}
In fact, set $K=e^{f}H_{0}$ and $f$ is defined by the Poisson equation
\begin{eqnarray*}
\triangle f= \frac{2}{r}
tr (\sqrt{-1}\Lambda_{\omega}(F_{H_{0}}+[\theta , \theta^{\ast H_{0}}
])-\lambda Id),
\end{eqnarray*}
by noting that $\int_{M }tr (\sqrt{-1}\Lambda_{\omega}(F_{H_{0}}+[\theta , \theta^{\ast H_{0}}])-\lambda Id) \frac{\omega ^{m}}{m!}=0$.

For simplicity, we denote:
\begin{eqnarray}
\Phi (H, \theta )=\sqrt{-1}\Lambda_{\omega}(F_{H}+[\theta , \theta^{\ast H}
])-\lambda Id_{E}.
\end{eqnarray}
It is easy to check that $(\sqrt{-1}\Lambda_{\omega}F_{H})^{\ast H}=\sqrt{-1}\Lambda_{\omega}F_{H}$ and $(\Phi (H, \theta ))^{\ast H}=\Phi (H, \theta )$. The following lemma is essentially proved by Simpson (\cite{Si}), we give a proof just for completeness.

\medskip

{\bf Lemma 4. } {\it Let $H(t)$ be a solution of the  heat flow (\ref{D1}) with initial metric $K$, then we have:
\begin{eqnarray}\label{F1}
(\frac{\partial }{\partial t}-\triangle )tr (\Phi (H, \theta ))=0
\end{eqnarray}
and
\begin{eqnarray}\label{F2}
(\frac{\partial }{\partial t}-\triangle )|\Phi (H, \theta )|_{H}^{2}=-4|D''_{ \theta} \Phi (H, \theta )|^{2}_{H},
\end{eqnarray}
where $D''_{ \theta}=\overline{\partial}_{E} +\theta $.}

\medskip

{\bf Proof. } Using the identities
\begin{eqnarray}
\begin{array}{lll}
&& \partial _{H}-\partial_{K} =h^{-1}\partial_{K}h ;\\
&& F_{H}-F_{K}=\overline{\partial }_{E} (h^{-1}\partial_{K} h) ;\\
&& \theta^{\ast H}=h^{-1} \theta^{\ast K} h ,\\
\end{array}
\end{eqnarray}
we have:
\begin{eqnarray}\label{F3}
\begin{array}{lll}
\frac{\partial }{\partial t }\Phi (H (t), \theta ) &=& \sqrt{-1} \Lambda_{\omega }\{ \overline{\partial }_{E} (-h^{-1} \frac{\partial h}{\partial t} h^{-1}\partial_{K} h +h^{-1}\partial _{K}(\frac{\partial h}{\partial t} ))\\
&& -[\theta ,  h^{-1}\frac{\partial h}{\partial t} h^{-1} \theta^{\ast K} h ]+ [\theta ,  h^{-1} \theta^{\ast K} \frac{\partial h}{\partial t}]\}\\
&=& \sqrt{-1} \Lambda_{\omega }\{ \overline{\partial }_{E} (\partial_{H}(h^{-1}\frac{\partial h}{\partial t}))+[\theta , [\theta^{\ast H} , h^{-1}\frac{\partial h}{\partial t}]\}.\\
\end{array}
\end{eqnarray}
The formula (\ref{F1}) can be deduced from (\ref{F3}) directly. On the other hand,
\begin{eqnarray}\label{F4}
\begin{array}{lll}
\triangle |\Phi |_{H}^{2}&=& -2\sqrt{-1}\Lambda_{\omega } \overline{\partial }\partial tr (\Phi H^{-1} \overline{\Phi }^{t} H)\\
&=& -2\sqrt{-1}\Lambda_{\omega } \overline{\partial } tr \{\partial \Phi H^{-1} \overline{\Phi }^{t} H- \Phi H^{-1} \partial H H^{-1} \overline{\Phi }^{t} H\\
& &+ \Phi H^{-1} \overline{\overline{\partial }\Phi }^{t} H + \Phi H^{-1} \overline{\Phi }^{t} H H^{-1}\partial H\\
&=& 2 Re <-2\sqrt{-1}\Lambda_{\omega } \overline{\partial }_{E}\partial_{E} \Phi , \Phi  >_{H} +2|\partial_{H}\Phi |_{H}^{2}+2|\overline{\partial}_{E}\Phi |_{H}^{2}\\
&& -2\sqrt{-1}\Lambda_{\omega } tr (\Phi H^{-1} \overline{[F_{H} , \Phi ]}^{t}H).
\end{array}
\end{eqnarray}
By (\ref{F3}) and (\ref{F4}), we have
\begin{eqnarray}\label{F4}
\begin{array}{lll}
&&(\triangle -\frac{\partial }{\partial t} ) |\Phi |_{H}^{2}=\triangle |\Phi |_{H}^{2} -2 Re <\frac{\partial }{\partial t}\Phi , \Phi  >_{H} \\
&=& 4 Re <\sqrt{-1}\lambda_{\omega }[\theta , [\theta^{\ast H} , \Phi ] , \Phi  >_{H}+ 2|\partial_{H}\Phi |_{H}^{2}+2|\overline{\partial}_{E}\Phi |_{H}^{2}\\
&=&2|\partial_{H}\Phi +[\theta^{\ast H} , \Phi ] |_{H}^{2}+2|\overline{\partial}_{E}\Phi + [\theta , \Phi ]|_{H}^{2}.\\
\end{array}
\end{eqnarray}
Since $\Phi (H, \theta )$ is self adjoint with respect to $H$, (\ref{F4}) implies (\ref{F2}).

\hfill $\Box$ \\

\medskip

Set $h=K^{-1}H=\exp (S)$, where $S\in End(E)$ and it is self adjoint with respect to $K $ or $H$.
By the initial condition and (\ref{F1}), we have $tr (\Phi (H(t), \theta ))=0$, then $\frac{\partial }{\partial t}det(h(t))=tr (h^{-1}\frac{\partial h}{\partial t})=0$.
So we have
\begin{eqnarray}
det (h(t))=1,
\end{eqnarray}
and
\begin{eqnarray}\label{1}
tr (S(t))=0.
\end{eqnarray}

By (\ref{F2}) and the maximum principle, we know that $\|\Phi (H(t), \theta )\|_{L^{\infty }}$ and $\|\Phi (H(t), \theta )\|_{L^{2 }}$ are monotonely decreasing.

For reader's convenience, we recall some notation. Let $K$ be a fixed Hermitian metric on bundle $E$, denote
\begin{eqnarray}\label{7.2}
S_{K}(E)=\{\eta \in \Omega^{0}(M, End(E))| \quad \eta ^{\ast K}=\eta \}.
\end{eqnarray}
Given $\rho \in C^{\infty }(R, R)$ and $\eta \in S_{K}(E)$. We define
$\rho (\eta )$ as follows. At each point $x$ on $M$, choose
an unitary basis $\{e_{i}\}_{i=1}^{r}$  with respect to metric
$K$, such that $\eta (e_{i})= \lambda_{i} e_{i}$. Set:
\begin{eqnarray}\label{7.3}
\rho (\eta ) (e_{i})= \rho (\lambda_{i}) e_{i} .
\end{eqnarray}

Given $\Psi \in C^{\infty } (R\times R , R)$, $\eta \in S_{K}(E)$,
$p\in \Omega^{0} (M, End(E))$. In a similar way, we define
$\Psi[\eta ](p)$ as follows. Let $\{e_{i}^{\ast }\}_{i=1}^{r}$ be the
dual basis for $\{e_{i}\}_{i=1}^{r}$, then $p\in \Omega^{0} (M,
End(E) )$ can be written
$$
p=\sum p_{ij} e_{i}^{\ast } \otimes e_{j}.
$$
Setting
\begin{eqnarray}\label{7.4}
\Psi [\eta ] (p)=\sum \Psi (\lambda _{i},
\lambda_{j})p_{ij}e_{i}^{\ast}\otimes e_{j}.
\end{eqnarray}
Let's recall  Donaldson's functional  defined on the space of Hermitian metrics on Higgs bundle $(E, \theta )$ ( Simpson \cite{Si}),
\begin{eqnarray}\label{7}
\mu (K, H) = \int_{M} tr (S \sqrt{-1}\Lambda_{\omega }F_{K, \theta })+ <\Psi (S) (D''_{\theta} S) , D''_{\theta} S>_{K}\frac{\omega^{n}}{n!},
\end{eqnarray}
where $\Psi (x, y)= (x-y)^{-2}(e^{y-x }-(y-x)-1)$, $\exp{S}=K^{-1}H$.
The Donaldson's functional has another equivalent form, i.e.
\begin{eqnarray}\label{F6}
\mu (K, H) =\int_{0}^{1} \int_{M}tr (\Phi (H(s), \theta ) H(s)^{-1}\frac{\partial H}{\partial s})\frac{\omega^{n}}{n!} ds,
\end{eqnarray}
where $H(s)$ is any path connecting metrics $K$ and $H$. The above integral is independent on   path, and we also have a formula for the derivative with respect
to $t$ of Donaldson¡¯s functional,
\begin{eqnarray}\label{F5}
\frac{d}{dt}\mu (K, H(t)) = \int_{M}tr (\Phi (H(s), \theta ) H(t)^{-1}\frac{\partial H}{\partial t})\frac{\omega^{n}}{n!},
\end{eqnarray}
for details see recent paper \cite{Ca}. By (\ref{F5}), we know that the heat flow (\ref{D1}) is the gradient flow of  Donaldson's functional.

\section{Proof of theorem 1 }
\setcounter{equation}{0}

Let $H(t)$ be the solution of the heat flow (\ref{D1}) with initial metric $K$, then we have
\begin{eqnarray}
\frac{d}{d t }\mu (K, H(t)) = -2\int_{M} |\Phi (H(t), \theta )|^{2}_{H(t)}\frac{\omega^{n}}{n!},
\end{eqnarray}
and
\begin{eqnarray}
-\mu (K, H(t)) =2\int_{0}^{t} \int_{M} |\Phi (H(s), \theta )|^{2}_{H(s)}\frac{\omega^{n}}{n!} ds .
\end{eqnarray}

\medskip

Case 1,   $\mu (K, H(t))\geq -C > -\infty $, i.e.
\begin{eqnarray}
\int_{0}^{+\infty} \int_{M} |\Phi (H(t), \theta )|^{2}_{H(t)}\frac{\omega^{n}}{n!} dt< \infty .
\end{eqnarray}
By the monotonicity of the integral, we have
\begin{eqnarray}\int_{M} |\Phi (H(t), \theta )|^{2}_{H(t)}\frac{\omega^{n}}{n!} \rightarrow 0
\end{eqnarray}
as $t\rightarrow +\infty$.

\medskip

Case 2,   $\mu (K, H(t)) \rightarrow -\infty $.

\medskip

By the definition of Donaldson's functional (\ref{7}), we have
\begin{eqnarray}\label{D3}
\mu (K, H(t)) \geq -C \|S(t)\|_{L^{\infty}}
\end{eqnarray}
where $\exp (S(t))=h(t)=K^{-1}H(t)$. From Lemma 4, we know that $\max _{M} |\sqrt{-1}\Lambda_{\omega } F_{H(t), \theta }|$ are uniformly bounded, so we have the following Simpson's estimate (p885 in \cite{Si})
\begin{eqnarray}\label{2}
\|S(t)\|_{L^{\infty}}\leq C_{1}\|S(t)\|_{L^{1}} + C_{2},
\end{eqnarray}
where constants $C_{1}$ and $C_{2}$ depend only on the curvature of initial metric $K$ and the geometry of $(M, \omega)$.
Then, (\ref{D3}) implies
\begin{eqnarray}
\|S(t)\|_{L^{1}}\rightarrow \infty
\end{eqnarray}
as $t \rightarrow \infty $.
By direct calculation, we have
\begin{eqnarray}
\begin{array}{lll}
&& \frac{\partial }{\partial t} \log (tr h(t) + tr h^{-1}(t))\\
&=& \frac{tr (h h^{-1}\frac{\partial h}{\partial t})- tr (h^{-1}\frac{\partial h }{\partial t}h^{-1})}{tr h + tr h^{-1}}\\
&\leq & |\Phi (H(t), \theta )|_{H(t)},\\
\end{array}
\end{eqnarray}
and
\begin{eqnarray}
\log (\frac{1}{2r}(tr h + tr h^{-1}))\leq |S|\leq r^{\frac{1}{2}}\log (tr h + tr h^{-1}),
\end{eqnarray}
where $r=rank (E)$.
By the above two inequalities, we have
\begin{eqnarray}
\begin{array}{lll}
&&r^{-\frac{1}{2}}\|S(t)\|_{L^{1}}- V\log 2r \\
&\leq & \int_{M} \log (tr h(t) + tr h^{-1}(t)) \frac{\omega^{n}}{n!} -V\log 2r
\\
&\leq & \int_{0}^{t} \|\Phi (H(s), \theta )\|_{L^{1}} ds\\
&\leq & t^{\frac{1}{2}}(\int_{0}^{t} \|\Phi (H(s) , \theta )\|_{L^{1}}^{2} ds )^{\frac{1}{2}}\\
&\leq & V^{\frac{1}{2}} t^{\frac{1}{2}}(\int_{0}^{t} \|\Phi (H(s) , \theta )\|_{L^{2}}^{2} ds )^{\frac{1}{2}}\\
&\leq &(V t)^{\frac{1}{2}}(- \mu (K , H(t)))^{\frac{1}{2}},
\end{array}
\end{eqnarray}
where $V$ is the volume of $(M, \omega )$. On the other hand, the monotonicity of $\|\Phi (H(s) , \theta )\|_{L^{2}}$ implies
\begin{eqnarray}
\begin{array}{lll}
&&t \|\Phi (H(t) , \theta )\|_{L^{2}}^{2}  \\
&\leq & \int_{0}^{t} \|\Phi (H(s) , \theta )\|_{L^{2}}^{2} ds \\
&= & -\mu (K , H(t))\\
\end{array}
\end{eqnarray}
 Combining the above two inequalities, we have
\begin{eqnarray}\label{key}
(r^{-\frac{1}{2}}\|S(t)\|_{L^{1}} -V\log 2r) \|\Phi (H(t) , \theta )\|_{L^{2}}\leq -V^{\frac{1}{2}}\mu (K, H(t)).
\end{eqnarray}

\medskip

{\bf Claim} {\it
Assume that the Higgs bundle $(E, \theta )$ is semi-stable, if $\mu (K, H(t))\rightarrow -\infty$, then we must have
\begin{eqnarray}
\lim _{t\rightarrow +\infty }\frac{-\mu (K, H(t))}{\|S(t)\|_{L^{1}}} =0.
\end{eqnarray}}

\medskip

By the above claim and (\ref{key}), we have $\|\Phi (H(t) , \theta )\|_{L^{2}}\rightarrow 0$ as $t\rightarrow +\infty$.
Combining the above two cases, we see that: if the Higgs bundle $(E, \theta )$ is semi-stable, then
\begin{eqnarray}
\|\Phi (H(t) , \theta )\|_{L^{2}}\rightarrow 0
\end{eqnarray}
as $t\rightarrow +\infty$.

Following Kobayashi's argument, let $u(x, t)=\int_{M} \chi (x, y, t-t_{0})|\Phi (H(t_{0}), \theta )|^{2} d \nu (y)$, where $\chi $ is the heat kernel. Using (\ref{F2}), we have
\begin{eqnarray}
(\frac{\partial }{\partial t} -\triangle ) (|\Phi (H(t), \theta )|_{H(t)}^{2}(x) -u(x, t))\leq 0.
\end{eqnarray}
By the maximum principle,
\begin{eqnarray*}
\max _{M} (|\Phi (H(t), \theta )|_{H(t)}^{2}(x) -u(x, t))\leq \max _{M} (|\Phi (H(t_{0}), \theta )|_{H(t_{0})}^{2}(x) -u(x, t_{0}))=0.
\end{eqnarray*}
Hence
\begin{eqnarray}
\begin{array}{lll}
\max _{M} |\Phi (H(t_{0}+1), \theta )|_{H}^{2}&\leq & \max _{M} u(x, t_{0}+1)\\
&=& \int_{M} \chi (x, y, 1)|\Phi (H(t_{0}), \theta )|_{H(t_{0})}^{2} d \nu (y)\\
&\leq & C\int_{M} |\Phi (H(t_{0}), \theta )|_{H(t_{0})}^{2} \frac{\omega ^{n}}{n!} \rightarrow 0,
\end{array}
\end{eqnarray}
as $t_{0} \rightarrow +\infty$. So there exists an approximate H-E metric structure on semi-stable Higgs bundle $(E, \theta )$.

To complete the proof of theorem 1, we only need to prove the above claim.

\medskip

{\bf Proof of the Claim. }
We will follow Simpson's argument (Proposition 5.3 in \cite{Si}) to show that if the estimate does not hold , there is a sub Higgs-sheaf contradicting semi-stability.

Suppose the required estimate does not hold. We can find a positive constant $C$ and a sequence $t_{i}\rightarrow +\infty$ such that
\begin{eqnarray}
\frac{-\mu (K, H(t_{i}))}{\|S(t_{i})\|_{L^{1}}}\geq C.
\end{eqnarray}

Set $u_{i}=l_{i}^{-1}S(t_{i})$, where $l_{i}=\|S(t_{i})\|_{L^{1}}\rightarrow +\infty$, then $\|u_{i}\|_{L^{1}}=1$.
By (\ref{1}) and (\ref{2}), we have $tr u_{i}=0$ and $\|u_{i}\|_{L^{\infty}}\leq C_{1}$. Simpson proved that:
 $u_{i}\rightarrow u_{\infty}$   weakly  in $L_{1}^{2}$; $\|u_{\infty}\|_{L^{1}}=1$, and the eigenvalues of $u_{\infty}$ are constant almost everywhere.

\medskip

Let $\lambda_{1} < \dots <\lambda _{l}$ denote the distinct eigenvalue of $u_{\infty}$.  Since $tr u_{\infty}=0$ and $\|u_{\infty}\|_{L^{1}}=1$, we must have $l\geq 2$.
For any $1\leq \alpha <l $, define function $P_{\alpha } : R\rightarrow R$ such that
\begin{eqnarray}
P_{\alpha }=\left \{\begin{array}{cll} 1, & x\leq \lambda_{\alpha }\\
0,
 & x\geq \lambda_{\alpha +1}\\
\end{array}\right.
\end{eqnarray}
Set $\pi_{\alpha }=P_{\alpha } (u_{\infty})$, Simpson (p887 in \cite{Si}) proved that:

\medskip

(1) $\pi_{\alpha } \in L_{1}^{2}$;

(2)  $\pi_{\alpha }^{2}=\pi_{\alpha }=\pi_{\alpha }^{\ast K}$;

(3) $(Id -\pi_{\alpha }) \bar{\partial }\pi_{\alpha } =0$;

(4) $(Id -\pi_{\alpha }) [\theta , \pi_{\alpha }]=0$.

\medskip

By Uhlenbeck and Yau's regularity statement of $L_{1}^{2}$-subbundle (\cite{UY}), $\pi_{\alpha }$ represent coherent torsion-free sub Higgs-sheaf $E_{\alpha }$ of $(E, \theta )$.
Set
\begin{eqnarray}
\nu =\lambda _{l} deg (E) -\sum_{\alpha =1} ^{l-1} (\lambda_{\alpha +1 } -\lambda_{\alpha }) deg (E_{\alpha}).
\end{eqnarray}
Since $u_{\infty }=\lambda _{l} Id  -\sum_{\alpha =1} ^{l-1} (\lambda_{\alpha +1} -\lambda_{\alpha })\pi_{\alpha }$ and $tr u_{\infty }=0$, we have
\begin{eqnarray}
\lambda _{l} rank (E) -\sum_{\alpha =1} ^{l-1}(\lambda_{\alpha +1} -\lambda_{\alpha }) rank (E_{\alpha })=0,
\end{eqnarray}
then
\begin{eqnarray}\label{3}
\nu =\sum_{\alpha =1} ^{l-1} (\lambda_{\alpha +1} -\lambda_{\alpha }) rank (E_{\alpha }) (\frac{deg (E)}{rank (E)}-\frac{deg (E_{\alpha })}{rank (E_{\alpha })}).
\end{eqnarray}

On the other hand, by the Gauss-Codazzi equation, we have
\begin{eqnarray}
deg (E_{\alpha })=\int_{M} tr (\pi_{\alpha } \sqrt{-1}\Lambda _{\omega }F_{K, \theta }) -|D''_{\theta} \pi_{\alpha }|_{K}^{2} \frac{\omega^{n}}{n!}.
\end{eqnarray}
Since $D''_{\theta} (\pi_{\alpha })=D''_{\theta}(P_{\alpha }(u_{\infty}))=dP_{\alpha } (u_{\infty})(D''_{\theta} u_{\infty
})$, where the function $df : R^{2} \rightarrow R$ defined by
\begin{eqnarray*}
df (x, y)= \frac{f(x)-f(y)}{x-y},
\end{eqnarray*}
and which is taken as $\frac{df }{dx}$ if $x=y$.
Then, we have
\begin{eqnarray}\label{6}
\begin{array}{lll}
\nu &=&\int_{M}tr ((\lambda _{l} Id  -\sum_{\alpha =1} ^{l-1}(\lambda_{\alpha +1} -\lambda_{\alpha }) \pi_{\alpha })\sqrt{-1}\Lambda _{\omega }F_{K, \theta })\\
&&+\sum_{\alpha =1} ^{l-1} (\lambda_{\alpha +1}-\lambda_{\alpha })|D''_{\theta} \pi_{\alpha }|^{2}\\
&=&\int_{M}tr (u_{\infty }\sqrt{-1}\Lambda _{\omega }F_{K, \theta })\\
&&+\langle \sum_{\alpha =1} ^{l-1} (\lambda_{\alpha +1}-\lambda_{\alpha })(dP_{\alpha })^{2}(u_{\infty }) (D''_{\theta} u_{\infty}) , D''_{\theta} u_{\infty}\rangle _{K}\\
\end{array}
\end{eqnarray}

If $\lambda_{k} \neq \lambda _{l}$, it easy to check that
\begin{eqnarray}\label{4}
\sum_{\alpha =1} ^{l-1} (\lambda_{\alpha +1}-\lambda_{\alpha })(dP_{\alpha })^{2}(\lambda_{k}, \lambda_{l})=|\lambda_{k}-\lambda_{l}|^{-1}
\end{eqnarray}
and
\begin{eqnarray}\label{5}
l_{i}\Psi (l_{i} \lambda_{k}, l_{i}\lambda_{l})\rightarrow =\left \{\begin{array}{cll}
(\lambda_{k}-\lambda_{l})^{-1},  & \lambda_{k} > \lambda_{l}\\
+\infty & \lambda_{k} \leq  \lambda_{l} .\\
\end{array}\right.
\end{eqnarray}
By (\ref{7}),  (\ref{4}), (\ref{5}) and (\ref{6}), we have
\begin{eqnarray}
\begin{array}{lll}
-C &\geq & \lim _{i\rightarrow \infty } \frac{\mu (K, H(t_{i}))}{\|S(t_{i})\|_{L^{1}}}\\
&=&\lim_{i\rightarrow \infty }\int_{M}tr (u_{i} \sqrt{-1}\wedge_{\omega }F_{K, \theta })+ \langle l_{i}\Psi (l_{i}u_{i}) (D''_{\theta} u_{i}) , D''_{\theta} u_{i}\rangle _{K}\\
&\geq &\int_{M}tr (u_{\infty }\sqrt{-1}\Lambda _{\omega }F_{K, \theta })\\
&&+\langle \sum_{\alpha =1} ^{l-1} (\lambda_{\alpha +1}-\lambda_{\alpha })(dP_{\alpha })^{2}(u_{\infty }) (D''_{\theta} u_{\infty}) , D''_{\theta} u_{\infty}\rangle _{K}\\
&=&\nu .\\
\end{array}
\end{eqnarray}
On the other hand,  (\ref{3}) and the semi-stability imply $\nu \geq 0$, so we get a contradiction.

\hfill $\Box$ \\

\hspace{0.3cm}

\end{document}